\documentclass[journal,twoside,web]{ieeecolor}
\usepackage{tabularx} 
\usepackage{diagbox}
\usepackage{graphicx}
\usepackage{mathrsfs}
\usepackage{generic}
\usepackage{cite}
\usepackage{amssymb,amsfonts}
\usepackage{algorithmic}
\usepackage{textcomp}
\def\BibTeX{{\rm B\kern-.05em{\sc i\kern-.025em b}\kern-.08em
    T\kern-.1667em\lower.7ex\hbox{E}\kern-.125emX}}
\usepackage[dvipsnames]{xcolor}
\usepackage{color, soul}
\sethlcolor{Green}

\usepackage{textcomp}
\def\BibTeX{{\rm B\kern-.05em{\sc i\kern-.025em b}\kern-.08em
    T\kern-.1667em\lower.7ex\hbox{E}\kern-.125emX}}
\usepackage{float}
\usepackage[tbtags]{amsmath}
\usepackage[tbtags]{mathtools}
\usepackage{accents}
\usepackage{multicol}
\usepackage{multirow}
\usepackage[ruled,vlined]{algorithm2e}
\usepackage{cite}
\usepackage{pgfplots}

\usepackage{array,multirow}
\usepackage{hhline}
\usepackage{arydshln}

\usepackage[hidelinks=true]{hyperref}

\newtheorem{theorem}{Theorem}

\newtheorem{lemma}{Lemma}
\newtheorem{example}{Example}
\newtheorem{assumption}{Assumption}

\newtheorem{problem}{Problem}

\newtheorem{remark}{Remark}

\makeatletter
\def\fnum@figure{\textcolor{subsectioncolor}{\sf Fig.~\thefigure}}
\def\fnum@table{\textcolor{subsectioncolor}{\sf TABLE~\thetable}}
\makeatother

\begin{document}
\title{Suboptimal Control of Unknown Second-order Nonlinear Systems with Guaranteed Global Convergence}
\author{Amir Shakouri, \IEEEmembership{Member, IEEE}, and M. Reza Emami, \IEEEmembership{Senior Member, IEEE}
\thanks{\textit{(Corresponding author: M. Reza Emami.)}}
\thanks{The authors are with the University of Toronto Institute for Aerospace Studies, Toronto, Canada (e-mail: \href{mailto:a_shakouri@outlook.com}{a$\_$shakouri@outlook.com}; \href{mailto:emami@utias.utoronto.ca}{emami@utias.utoronto.ca}).}}

\maketitle

\pagestyle{empty}
\thispagestyle{empty}

\begin{abstract}
A suboptimal active disturbance rejection controller (S-ADRC) is proposed for second-order systems with unknown time-varying nonlinear dynamics. The output-feedback controller guarantees a global convergence to the vicinity of an optimal solution by means of dynamic control gains, based on the estimated main and extended state variables obtained through a high-gain observer. Three numerical examples compare the performance of the proposed control scheme applied to linear and nonlinear systems with that of a fixed-gain conventional ADRC as well as several model-based optimal and suboptimal controllers. 
\end{abstract}

\section{Introduction}
\label{sec:I}

\IEEEPARstart{T}{he} control of nonlinear systems with unknown dynamics becomes more challenging when there is a high priority for energy consumption while little or no rich experimental data is available for the system. Numerous critical engineering systems fall into this category. For instance, robot manipulators in an assembly plant may continuously need to handle variable unknown payloads with minimal energy consumption to reduce both the overall production costs and the $\mathrm{CO_2}$ emissions in the long-term \cite{meike2013energy}. Another example is the deployment of a (small) spacecraft for removing a piece of debris from near-Earth environment \cite{hakima2018assessment} or redirecting the natural orbit of a near-Earth asteroid \cite{bazzocchi2016comparative}, by attaching the spacecraft to the object to steer it away, where the effect of the unknown object on the system dynamics is highly uncertain, while the available energy onboard the spacecraft is quite limited. 

The optimal and suboptimal controllers are mainly based on a fairly accurate system model \cite{kirk2004optimal,mracek1998control,ccimen2008state,erdem1999globally,chen2015global,lin2018control}. In the absence of such a model, the existing model-free controllers need to access rich data for an online or offline learning phase \cite{khan2012reinforcement,kiumarsi2017optimal,hou2011data,werbos1989neural,verhaegen2007filtering} or for implementing a direct data-driven approach \cite{de2019formulas,van2020noisy,nortmann2020data,guo2021data,dai2021nonlinear}. If such an experimental database is not available, to compensate for nonlinear unmodeled dynamics, high-gain methods have been proposed in the literature
\cite{zheng2007stability,shakouri2021framework}, which are mostly quite energy demanding. Without an accurate system model, if the (time-varying) nonlinear dynamics of the system plays a significant role, the state-of-the-art control schemes hardly guarantee the global convergence of the system to the vicinity of an energy optimal solution. 

The implementation of active disturbance rejection controller (ADRC) \cite{zheng2007stability} along with the extended state observer (ESO) \cite{guo2011convergence} has been a popular approach to constructing an output-feedback scheme capable of regulating unknown nonlinear systems. However, the conventional ADRC may not be a proper solution for many practical systems, as it demands a high-energy effort emanating from its feedback linearization nature. 

In this letter, we make use of the ESO to construct a suboptimal output-feedback control law, called \textit{suboptimal active disturbance rejection controller} (S-ADRC), for double-integrator systems with unknown Lipschitz dynamics. It is shown that the conventional ADRC is a special type of the proposed scheme. The S-ADRC is suboptimal in the sense that its closed-loop response globally converges to the vicinity of an optimal solution. Therefore, one may argue that the proposed S-ADRC has the advantages of both the conventional ADRC and the state-dependent Riccati equation (SDRE) controller, as the former globally regulates a nonlinear system in the presence of unmodeled dynamics \cite{zheng2007stability}, and the latter is a model-based suboptimal controller \cite{ccimen2008state}.

\section{Preliminaries}
\label{sec:II}

\subsection{Notations}
\label{sec:II-1}

Let $\mathbb{R}^n$ and $\mathbb{R}^{m,n}$ denote the space of $n$-dimensional vectors and $m \times n$ real matrices, respectively. The $i$th entry of vector $r$ and the $ij$th entry of matrix $M$ are referred to by $r_i$ and $M_{ij}$, respectively. The inverse and transpose of matrix $M$ are denoted by $M^{-1}$ (if the inverse exists) and $M^{\top}$, respectively. The symbol $\|\cdot\|$ denotes the 2-norm for vectors, and we use $|\cdot|$ to show the absolute value for scalars. A vector-valued function $v$ is said to be in the order of $\varepsilon$, denoted by $\|v\|=\mathrm{O}(\varepsilon)$, when there exists a real constant $\sigma>0$ such that $\|v\|\leq\sigma\varepsilon$. We use $\min\{a,b\}$ to denote the smallest value between $a$ and $b$. The class of $k$-times differentiable functions is denoted by $\mathcal{C}^k$.

\subsection{Problem statement}

Consider a second-order system as follows: 
\begin{equation}
\label{eq:1}
\left\{\begin{array}{lcl}
\dot{x}_1&=&x_2, \\
\dot{x}_2&=&f(x_1,x_2,u,t)+gu, \\
y&=&x_1,
\end{array}\right.
\end{equation}
where $x_1,x_2\in\mathbb{R}$ are the state variables, and $y\in\mathbb{R}$ is the output of the system. The control input $u\in\mathcal{C}^1$ depends on the output $y$ through an observer. The input gain $g$ and the system dynamics $f$ satisfy the following condition. 

\begin{assumption}
\label{ass:1}
The input gain $g$ is a known nonzero real constant. Function $f(x_1,x_2,u,t)$ is unknown with a bounded rate, i.e., there exists a $\varphi>0$ such that $|\dot{f}(x_1,x_2,u,t)|\leq\varphi$ for all $x_1,x_2,u\in\mathbb{R}$, and $t\in[0,\infty)$. Moreover, for any $\sigma>0$ there exists a $\delta>0$ such that if $|x_1|+|x_2|<\sigma$, then $|f(x_1,x_2,u,t)|<\delta$ for all $t\in[0,\infty)$ and $u\in\mathbb{R}$.
\end{assumption}

Note that the results of this letter can be extended to multi-variable cases of system \eqref{eq:1}, where $x_1$, $x_2$, $f$, and $y$ are vectors, and $g$ is replaced by a known diagonal matrix. However, for simplicity and brevity, we build our analysis on system \eqref{eq:1}. See Remark \ref{rem:3} and Example \ref{ex:3} for more detail. 

Consider the following infinite-horizon quadratic cost function corresponding to system \eqref{eq:1}:
\begin{equation}
\label{eq:III-8}
J(x_1,x_2,u)=\frac{1}{2}\int_{0}^\infty(q_1x_1^2+q_2x_2^2+ru^2)\mathrm{d}t,
\end{equation}
where $q_1,q_2,r>0$ are constant weight factors. It is well-known that the closed-loop solution of system \eqref{eq:1} under an optimal control law corresponding to the cost function \eqref{eq:III-8}, called the \textit{optimal solution}, satisfies the following conditions \cite{kirk2004optimal}:
\begin{subequations}
\label{eq:III-9}
\begin{align}
\label{eq:III-9a}
\frac{\partial}{\partial u}H(x_1,x_2,u,\lambda)&=0,\\
\label{eq:III-9b}
\dot{\lambda}_i+\frac{\partial}{\partial x_i}H(x_1,x_2,u,\lambda)&=0,\ i=1,2,
\end{align}
\end{subequations}
where $\lambda_1$ and $\lambda_2$ are the adjoint variables, and $H$ is the Hamiltonian corresponding to system \eqref{eq:1} and cost function \eqref{eq:III-8}, defined as follows:
\begin{equation}
\label{eq:III-10}
\begin{split}
H(x_1,x_2,u,\lambda)=&\frac{1}{2}(q_1x_1^2+q_2x_2^2+ru^2)+\lambda_1x_2\\
&+\lambda_2(f(x_1,x_2,u,t)+gu).
\end{split}
\end{equation}

The control problem can be stated as follows:

\begin{problem}
\label{prob:1}
Having  Assumption \ref{ass:1}, find an output-feedback control law under which the closed-loop response of the unknown system \eqref{eq:1} globally converges to and remains in the vicinity of an optimal solution.
\end{problem}

\section{Results}    
\label{sec:III}

Let us define an extended state variable representing the unknown dynamics of the system as $x_3=f(x_1,x_2,u,t)$. We use the solution of a linear extended state observer (LESO) \cite{guo2011convergence} to estimate $x_1$, $x_2$ and the extended state variable $x_3$ for system \eqref{eq:1}, based on the output $y$. (Nonlinear extended state observers \cite{guo2011convergence} or other disturbance observers with the same objective \cite{chen2015disturbance} can also be used for this purpose.) The estimates are then used in the formulation of the proposed controller. Accordingly, we define the total estimated state as $\hat{x}=[\hat{x}_1,\ \hat{x}_2,\ \hat{x}_3]^{\top}$, and solve the following system starting from some initial condition $\hat{x}(0)$:
\begin{equation}
\label{eq:III-1}
\left\{
\begin{array}{lcl}
\dot{\hat{x}}_1&=&\hat{x}_2+(\alpha_1/\varepsilon)(y-\hat{x}_1), \\
\dot{\hat{x}}_2&=&\hat{x}_3+(\alpha_2/\varepsilon^2)(y-\hat{x}_1)+gu, \\
\dot{\hat{x}}_3&=&(\alpha_3/\varepsilon^3)(y-\hat{x}_1),
\end{array}\right.
\end{equation}
where $0<\varepsilon<1$ is constant and $\alpha_i$, $i=1,2,3$ are selected such that the following matrix $E$ is Hurwitz:
\begin{equation}
\label{eq:III-2}
E=\begin{bmatrix}
-\alpha_1 & 1 & 0 \\
-\alpha_2 & 0 & 1 \\
-\alpha_3 & 0 & 0
\end{bmatrix}.
\end{equation}

The proposed S-ADRC for system \eqref{eq:1} can be expressed as follows:
\begin{equation}
\label{eq:III-3}
u=-\frac{1}{g}(p\hat{x}_1+d(\hat{x})\hat{x}_2+\hat{x}_3),
\end{equation}
where $p>0$ is constant, and $d(\hat{x})$ is a state-dependent function obtained as
\begin{equation}
\label{eq:III-4}
d(\hat{x})=\mathrm{min}\left\{\delta,\sqrt{\left(\frac{\hat{x}_3-\rho\hat{x}_1}{\hat{x}_2}\right)^2+\beta}\right\},
\end{equation}
where $\rho$ and $\beta$ are constant control gains with the following conditions:
\begin{subequations}
\begin{align}
\label{eq:III-5a}
0<\rho&<p,\\
\label{eq:III-5b}
2(p+\rho)&<\beta.
\end{align}
\end{subequations}

\begin{remark}
\label{rem:1}
Note that if $\beta=\delta^2$, then the derivative gain of the S-ADRC is constant $d=\delta$. In this case, the S-ADRC reduces to a conventional ADRC with constant gains \cite{zheng2007stability}. 
\end{remark}

\subsection{S-ADRC Convergence}

The following theorem presents the convergence behavior of the proposed S-ADRC:

\begin{theorem}
\label{th:stablity}
Suppose that Assumption \ref{ass:1} holds. Let parameters $\alpha_i$, $i=1,2,3$, of the LESO \eqref{eq:III-1} be selected such that matrix $E$ expressed in \eqref{eq:III-2} is Hurwitz. Then, the unknown system \eqref{eq:1} under the S-ADRC scheme, consisting of observer \eqref{eq:III-1} and controller \eqref{eq:III-3}, satisfies the following attractivity and boundedness properties for any bounded initial condition $\hat{x}(0)$ (and $x(0)$):
\begin{enumerate}
\item There exists a $\sigma>0$ such that $\|x(t)\|<\sigma$ for all $t\geq{0}$.
\item $x\rightarrow0$ as $\varepsilon \rightarrow0$ and $t\rightarrow\infty$.
\item There exists a $t_0\geq0$ such that $\|x-\hat{x}\|=\mathrm{O}(\varepsilon)$, $\forall t\geq t_0$.
\end{enumerate}\hfill$\bullet$
\end{theorem}

According to Theorem \ref{th:stablity}, the convergence characteristics of the S-ADRC are similar to those of the conventional ADRC with constant gains \cite[Thm. 4]{zheng2007stability}, \cite[Thm. 2.1]{guo2011convergence}. 

\subsubsection*{Proof of Theorem \ref{th:stablity}} First, we need the following lemma:
\begin{lemma}
\label{lem:b}
Consider a perturbed stable linear time-invariant (LTI) system $\dot{z}(t)=\theta_1 Z z(t)+\theta_2\zeta(t)$ where $Z$ is a Hurwitz matrix, $\theta_1,\theta_2>0$, and $\zeta(t)$ is a bounded perturbation such that $\|\zeta(t)\|<\bar{\zeta}$ for all $t\in[0,\infty)$. Then, there exist $\lambda<0$ and $c>0$ such that the solution of the system $z(t)$ is bounded as follows:
\begin{equation}
\label{eq:lem:b}
\|z(t)\|\leq c\left(\exp(\theta_1\lambda t)\|z(0)\|-\frac{\theta_2\bar{\zeta}}{\theta_1\lambda}\right).
\end{equation}\hfill$\bullet$
\end{lemma}

\begin{proof}
According to the variation of constants formula \cite[p. 74]{coddington1955theory}, the solution of $\dot{z}(t)=\theta_1Z z(t)+\theta_2\zeta(t)$ can be written as
\begin{equation}
\label{eq:prof-3}
z(t)=\exp(\theta_1Zt)z(0)+\int_0^{t}\exp(\theta_1Z(t-s))\theta_2\zeta(s)\mathrm{d}s,
\end{equation}
and the following inequality can be directly concluded:
\begin{equation}
\label{eq:prof-4}
\|z(t)\|\leq\|\exp(\theta_1Zt)z(0)\|+\theta_2\int_0^{t}\|\exp(\theta_1Z(t-s))\zeta(s)\|\mathrm{d}s.
\end{equation}
It is known that for every scalar $s$ and vector $v$ there exist $c>0$ and $\lambda<0$ such that $\|\exp(\theta_1Z(t-s))v\|\leq c\exp(\theta_1\lambda(t-s))\|v\|$ (see \cite[Thm. 2.34]{chicone2006ordinary}). Therefore, since we have $\|\zeta(t)\|\leq\bar{\zeta}$, inequality \eqref{eq:prof-4} reduces to
\begin{equation}
\label{eq:prof-5}
\begin{split}
\|z(t)\|&\leq c\exp(\theta_1\lambda t)\|z(0)\|+c\theta_2\bar{\zeta} \int_0^{t}\exp(\theta_1\lambda(t-s))\mathrm{d}s\\
&\leq c\exp(\theta_1\lambda t)\|z(0)\|-\frac{c\theta_2\bar{\zeta}}{\lambda}\left(1-\exp(\theta_1\lambda t)\right)\\
&\leq c\left(\exp(\theta_1\lambda t)\|z(0)\|-\frac{\theta_2\bar{\zeta}}{\theta_1\lambda}\right).
\end{split}
\end{equation}
\end{proof} 

We define the estimation error of the LESO as $\tilde{x}=x-\hat{x}$ and the normalized estimation error as $e_1=\tilde{x}_1$, $e_2=\varepsilon\tilde{x}_2$, and $e_3=\varepsilon^2\tilde{x}_3$. The normalized error vector $e=[e_1,\ e_2,\ e_3]^{\top}$ obeys the following dynamics:
\begin{equation}
\label{eq:prof-2}
\dot{e}=\frac{1}{\varepsilon}Ee+\varepsilon^2h(t),
\end{equation}
where $h(t)=[0,\ 0,\ \dot{f}(x_1,x_2,u,t)]^{\top}$. According to Assumption \ref{ass:1} we have $\|h(x,t,u)\|\leq\varphi$. Thus, Lemma \ref{lem:b} can be used to verify the following inequality:
\begin{equation}
\label{eq:prof-5}
\begin{split}
\|e(t)\|\leq c\left(\exp(\lambda t/\varepsilon)\|e(0)\|-\frac{\varphi\varepsilon^3}{\lambda}\right).
\end{split}
\end{equation}
Since $\exp(\lambda t/\varepsilon)\|e(0)\|$ is a strictly decreasing function of time, for every $\sigma>0$ there exists a $t_0>0$ such that $\exp(\lambda t/\varepsilon)\|e(0)\|\leq\sigma\varepsilon^3$ for all $t\geq t_0$. Therefore, there exists a $t_0>0$ such that for all $t\geq t_0$ we have:
\begin{equation}
\label{eq:prof-5p}
\|e(t)\|\leq c\left(\sigma-\varphi/\lambda\right)\varepsilon^3.
\end{equation}
Since $\tilde{x}_i=e_i/\varepsilon^{i-1}$, $i=1,2,3$, given that $\varepsilon<1$ and using \eqref{eq:prof-5p}, the following inequality holds:
\begin{equation}
\label{eq:prof-6}
\begin{split}
\|\tilde{x}(t)\|\leq\frac{\|e(t)\|}{\varepsilon^2}\leq c\left(\sigma-\varphi/\lambda\right)\varepsilon,\ \forall t>t_0,
\end{split}
\end{equation}
which proves the third property of Theorem \ref{th:stablity}.

The closed-loop dynamics, constructed by system \eqref{eq:1}, the LESO \eqref{eq:III-1}, and the proposed S-ADRC \eqref{eq:III-3}, can be represented as:
\begin{equation}
\label{th:stablity-2}
\left\{
\begin{array}{lcl}
\dot{x}_1&=&x_{2}, \\
\dot{x}_{2}&=&-px_1-d(x)x_{2}+\underbrace{\tilde{x}_3+p\tilde{x}_1+\tilde{d}}_{\upsilon},
\end{array}\right.
\end{equation}
where $\tilde{d}=d(x)x_2-d(\hat{x})\hat{x}_{2}$, and $|\tilde{d}|$ is bounded, as $|x_2-\hat{x}_{2}|$ is bounded (from \eqref{eq:prof-6}) and $\sqrt{\beta}\leq d(\cdot) \leq \delta$. Therefore, $\upsilon$ is bounded, and from \eqref{eq:prof-6} we have $\upsilon\rightarrow0$ as $\varepsilon\rightarrow0$. 
Consider the following candidate Lyapunov function: 
\begin{equation}
\label{th:stablity-3}
V=\frac{1}{2}(px_1^2+x_2^2),
\end{equation}
and verify that
\begin{equation}
\label{th:stablity-4}
\begin{split}
\dot{V}&=px_1\dot{x}_1+x_2\dot{x}_2\\
&=px_1x_2+x_2(-px_1-d(x)x_2+\upsilon)\\
&=-d(x)x_2^2+\upsilon x_2\leq-\sqrt{\beta}x_2^2+|\upsilon||x_2|.
\end{split}
\end{equation}
Given \eqref{th:stablity-3} and \eqref{th:stablity-4}, one can verify that $|x_2|$ is upper-bounded, because if it keeps on increasing, then $\dot{V}$ becomes negative at some point, indicating that $V$ is strictly decreasing, which in turn shows that $|x_2|$ will be decreasing, as well. Therefore, $|x_2|$ must be limited to an upper bound.

System \eqref{th:stablity-2} can be written in the following form:
\begin{equation}
\label{th:stablity-5}
\begin{bmatrix}
\dot{x}_1 \\
\dot{x}_2
\end{bmatrix}=\begin{bmatrix}
0 & 1 \\
-p & -\sqrt{\beta}
\end{bmatrix}\begin{bmatrix}
x_1 \\
x_2
\end{bmatrix}+\begin{bmatrix}
0 \\
\omega
\end{bmatrix},
\end{equation}
where $\omega=\upsilon+(\sqrt{\beta}-d(x))x_2$ is bounded as $\upsilon$, $d(x)$, and $x_2$ are bounded. Given Lemma \ref{lem:b}, we conclude that since  system \eqref{th:stablity-5} is governed by a stable LTI system perturbed by a bounded additive term, the state vector is bounded, which proves the first property of Theorem \ref{th:stablity}. For the second property of Theorem \ref{th:stablity}, recall that $\upsilon\rightarrow0$ as $\varepsilon\rightarrow0$, which implies that in \eqref{th:stablity-4} $\dot{V}$ becomes strictly negative, and hence $x\rightarrow0$ as $\varepsilon\rightarrow0$ and $t\rightarrow\infty$. \hfill$\QED$

\subsection{S-ADRC Suboptimality}

In the context of SDRE control, a known system is called suboptimal if \eqref{eq:III-9a} always holds and \eqref{eq:III-9b} is eventually satisfied. In the case of this analysis, as we deal with an unknown system, we call a system \textit{suboptimal} if \eqref{eq:III-9a} always approximately holds and \eqref{eq:III-9b} is eventually approximately satisfied.

The following theorem states the suboptimality of the S-ADRC:

\begin{theorem}
\label{th:optimality}
Suppose that Assumption \ref{ass:1} holds. Let parameters $\alpha_i$, $i=1,2,3$, of the LESO \eqref{eq:III-1} be selected such that matrix $E$ expressed in \eqref{eq:III-2} is Hurwitz. Let the weight factors of the cost function $J(x_1,x_2,u)$, described in \eqref{eq:III-8}, satisfy the following conditions:
\begin{subequations}
\begin{align}
\label{eq:th:optimality-1a}
\frac{q_1}{r}&=\frac{p^2-\rho^2}{g^2},\\
\label{eq:th:optimality-1b}
\frac{q_2}{r}&=\frac{\beta-2\left(p+\rho\right)}{g^2},
\end{align}
\end{subequations}
which are adjustable to any desired positive values by a proper selection of $p$, $\rho$, and $\beta$. Also, let the set $\mathcal{X}\subset\mathbb{R}^3$ be defined as follows:
\begin{equation}
\label{eq:III-7}
\begin{split}
\mathcal{X}=\left\{\hat{x}\in\mathbb{R}^3:\sqrt{\left(\frac{\hat{x}_3-\rho\hat{x}_1}{\hat{x}_2}\right)^2+\beta}\leq\delta\right\}.
\end{split}
\end{equation}
Then, as long as $\hat{x}\in\mathcal{X}$, the response of the unknown system \eqref{eq:1} under the S-ADRC scheme, consisting of the observer \eqref{eq:III-1} and controller \eqref{eq:III-3}, is suboptimal in the sense that it satisfies the following two conditions:
\begin{subequations}
\begin{align}
\label{eq:th:optimality-2a}
\exists t_0>0:\ \left|\frac{\partial}{\partial u}H(x_1,x_2,u(\hat{x}),\lambda)\right|&\leq \mathrm{O}(\varepsilon),\ \forall t\geq t_0,\\
\label{eq:th:optimality-2b}
\lim_{t\rightarrow\infty}\left|\dot{\lambda}_1+\frac{\partial}{\partial x_1}H(x_1,x_2,u(\hat{x}),\lambda)\right|&\leq \mathrm{O}(\varepsilon),\ i=1,2.
\end{align}
\end{subequations}
\hfill$\bullet$
\end{theorem}

\subsubsection*{Proof of Theorem \ref{th:optimality}} Note that according to Theorem 1 and Assumption 1, as $x$ is bounded under the S-ADRC, $f(x_1,x_2,u,t)$ should be bounded. Also, since based on (17) the state error $\tilde{x}$ is bounded, $\hat{x}$ is also bounded, which by virtue of (7) results in $u$ being bounded. Therefore, the cost function (2) remains bounded, which reduces to the existence of an optimal solution \cite[\S 3.1]{ccimen2008state}. 

Observe that the state variables of system \eqref{eq:1} can be parameterized in terms of state-dependent coefficient (SDC) matrix $A(x)$ as:
\begin{equation}
\label{eq:th:opt-4}
\begin{bmatrix}
\dot{x}_1 \\
\dot{x}_2 
\end{bmatrix}=A(x)\begin{bmatrix}
x_1 \\
x_2
\end{bmatrix}+Bu,
\end{equation}
where
\begin{equation}
\label{eq:th:opt-5}
A(x)=\begin{bmatrix}
0 & 1  \\
\rho & \frac{x_3-\rho x_1}{x_2} 
\end{bmatrix},\
B=\begin{bmatrix}
0 \\
g
\end{bmatrix}.
\end{equation}
Define $Q=\mathrm{diag}(q_1,q_2)$, and let $P(x)$ be the solution of the following SDRE:
\begin{equation}
\label{eq:th:opt-6}
\begin{split}
A^{\top}(x)P(x)+P(x)A(x)-\frac{1}{r}P(x)BB^{\top}P(x)+Q=0,
\end{split}
\end{equation}
which is a symmetric positive-definite matrix. The solution of the SDRE \eqref{eq:th:opt-6} can be expressed in the following entry-wise form \cite{chen2015global}:
\begin{subequations}
\begin{align}
\label{eq:th:opt-7a}
P_{12}=&\frac{r\rho}{g^2}+\sqrt{\frac{r^2\rho^2}{g^4}+\frac{rq_1}{g^2}},\\
\label{eq:th:opt-7b}
P_{22}=&\frac{r(x_3-\rho x_1)}{g^2x_2}+\sqrt{\frac{r^2}{g^4}\left(\frac{x_3-\rho x_1}{x_2}\right)^2+\frac{rq_2}{g^2}+\frac{2rP_{12}}{g^2}}.
\end{align}
\end{subequations}
Using conditions \eqref{eq:th:optimality-1a} and \eqref{eq:th:optimality-1b}, equations \eqref{eq:th:opt-7a} and \eqref{eq:th:opt-7b} can be reduced to:
\begin{subequations}
\begin{align}
\label{eq:th:opt-7ap}
P_{12}&=\frac{r}{g^2}(\rho+p),\\
\label{eq:th:opt-7bp}
P_{22}&=\frac{r}{g^2}\left[\frac{x_3-\rho x_1}{x_2}+\sqrt{\left(\frac{x_3-\rho x_1}{x_2}\right)^2+\beta}\right].
\end{align}
\end{subequations}

Verify that when $x\in\mathcal{X}$ and $\hat{x}=x$, the proposed S-ADRC in \eqref{eq:III-3} reduces to the following form in terms of $P(x)$ entries:
\begin{equation}
\label{eq:th:opt-8}
u(x)=-\frac{g}{r}(P_{12}x_1+P_{22}(x)x_2),
\end{equation}
which is essentially an SDRE controller (see \cite[Eqs. (8)--(14)]{chen2015global}). 

We know that the adjoint variables of the optimal solution satisfy the following condition (see \cite[Eq. (54)]{mracek1998control}, \cite[Eqs. (16), (22)]{ccimen2008state}):
\begin{equation}
\label{eq:th:opt-8p}
\begin{bmatrix}
\lambda_1 \\
\lambda_2
\end{bmatrix}
=P(x)\begin{bmatrix}
x_1\\
x_2
\end{bmatrix}.  
\end{equation}

\begin{lemma}[see {\cite[Thm. 5]{ccimen2008state}}, {\cite[Thm. 2]{mracek1998control}}]
\label{lem:1}
For the nominal case, system \eqref{eq:1} under SDRE controller \eqref{eq:th:opt-8} satisfies the first property of optimality described by \eqref{eq:III-9a}, that is, $\partial H(x_1,x_2,u,\lambda)/\partial u=0$.\hfill$\bullet$
\end{lemma}

The situation for the S-ADRC system is slightly different from the nominal case stated in Lemma \ref{lem:1}, because the controller is characterized as:
\begin{equation}
\label{eq:th:opt-8pp}
u(\hat{x})=-\frac{g}{r}(P_{12}\hat{x}_1+P_{22}(\hat{x}) \hat{x}_2),
\end{equation}
and provided that $P(x)$ is Lipschitz continuous, we can substitute \eqref{eq:th:opt-8p} and \eqref{eq:th:opt-8pp} in \eqref{eq:III-9a} to obtain:
\begin{equation}
\label{eq:th:opt-9}
\begin{split}
\left|\frac{\partial}{\partial u}H(x_1,x_2,u(\hat{x}),\lambda)\right|=&\big|P_{12}(x_1-\hat{x}_1)+P_{22}(x)x_2\\
&-P_{22}(\hat{x})\hat{x}_2\big|\leq \mathrm{O}(\varepsilon),
\end{split}
\end{equation}
which shows the validity of \eqref{eq:th:optimality-2a}. 

For proving the second property, consider the following lemma for the nominal system:
\begin{lemma}[see {\cite[Thm. 5]{ccimen2008state}}, {\cite[Thm. 4]{mracek1998control}}]
\label{lem:2}
For the nominal case, the solution of system \eqref{eq:1} under SDRE controller \eqref{eq:th:opt-8}, asymptotically (at a quadratic rate) approaches \eqref{eq:III-9b}. More precisely, there exists a constant positive-definite matrix $M\in\mathbb{R}^{6\times 6}$ by which the following inequality holds:
\begin{equation}
\label{eq:th:opt-10}
\begin{split}
\left|\dot{\lambda}_i+\frac{\partial H(x_1,x_2,u(x),\lambda)}{\partial x_i}\right|\leq x^{\top}Mx,\ i=1,2.
\end{split}
\end{equation}\hfill$\bullet$
\end{lemma}

To extend Lemma \ref{lem:2} to the system with controller $u(\hat{x})$, consider the following inequalities: 
\begin{equation}
\label{eq:th:opt-11}
\begin{split}
\left|\frac{\partial H(x_1,x_2,u(\hat{x}),\lambda)}{\partial x_i}-\frac{\partial H(x_1,x_2,u(x),\lambda)}{\partial x_i}\right|&\leq \mathrm{O}(\|\hat{x}-x\|)\\&\leq \mathrm{O}(\varepsilon).
\end{split}
\end{equation}
One can verify that the first inequality holds according to the Lipschitz continuity of all terms, and the second inequality holds according to the estimation error bound of the LESO obtained in \eqref{eq:prof-6}. Therefore, since based on \eqref{eq:th:opt-10} one of the terms in the left-hand side of \eqref{eq:th:opt-11} asymptotically approaches $-\dot{\lambda}_i$ at a quadratic rate, we can conclude that $|\dot{\lambda}_i+\partial H(x_1,x_2,u(\hat{x}),\lambda)/\partial x_i|$ is bounded by $\mathrm{O}(\varepsilon)$ as time approaches infinity, which is obviously reduced to the second property of suboptimality \eqref{eq:th:optimality-2b}.\hfill$\QED$

\begin{remark}
\label{rem:3}
The proposed S-ADRC inherits the suboptimality feature from its similarities to the SDRE controller \cite{ccimen2008state,mracek1998control}. The main advantages of the S-ADRC over the SDRE controller are that a) it does not use any model of the system, and b) it is globally convergent. Although SDRE controllers can be designed to be globally asymptotically stable for single-variable second-order systems \cite{erdem1999globally,chen2015global,lin2018control}, this property cannot be extended easily to multi-variable cases. However, the proposed S-ADRC can be implemented in second-order multi-variable systems, such as many mechanical systems (see Example \ref{ex:3}). 
\end{remark}

\section{Examples}    
\label{sec:VI}

In this section, we demonstrate some examples of the proposed S-ADRC for unknown linear time-varying (LTV) and nonlinear systems\footnote{MATLAB\textsuperscript{\tiny\textregistered} codes and Simulink\textsuperscript{\tiny\textregistered} models for the examples of this letter can be found in \href{https://github.com/a-shakouri/s-adrc}{https://github.com/a-shakouri/s-adrc}}. To assess the optimality of the proposed scheme, we compare the S-ADRC response with the optimal/suboptimal controllers obtained via a full knowledge of the system model and state variables. In addition, we compare the results with the most efficient conventional constant-gain ADRC. 

\begin{figure}[!h]
\centering\includegraphics[width=1\linewidth]{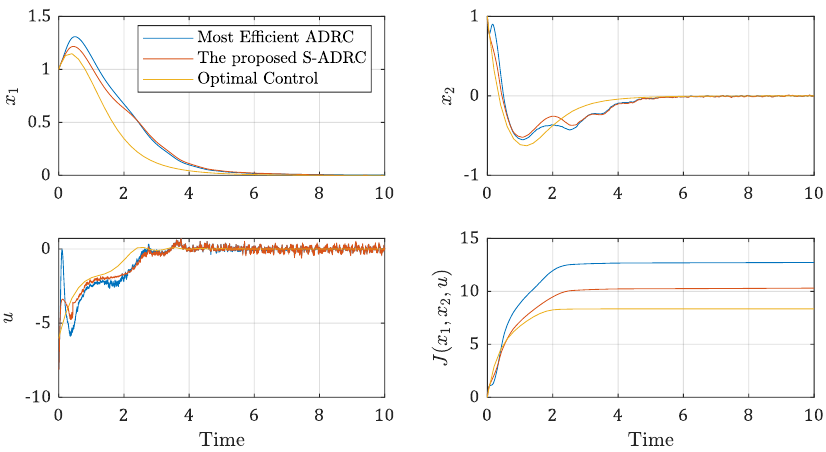}
\caption{Simulation results for Example \ref{ex:1}.}
\label{fig:1}
\end{figure}

\begin{example}
\label{ex:1}
Consider a second-order LTV system in the form of \eqref{eq:1} with $g=1$ and $f(x_1,x_2,t)=\exp(\sin(t))x_1+\cos(t^2)x_2$, which has an unstable open-loop dynamics. We select the LESO parameters as $\alpha_1=\alpha_2=3$, $\alpha_3=1$, and $\varepsilon=0.05$ in order to make matrix $E$ Hurwitz with all its eigenvalues equal to $-1$. The control parameters are considered as $p=2$, $\rho=1$, and  $\beta=9$ in order to obtain equal weight factors of $q_1/r=q_2/r=3$. Also, a large value is selected for parameter $\delta=10^6$ to enlarge the suboptimality regime. We corrupt the measurement output as $y=x_1+\nu$, where the noise $\nu$ is a zero-mean normally distributed random variable with standard deviation of $0.01$. The initial condition of the system is $x(0)=[1,\ 1]^{\top}$, and initial estimated state for the observer is set to be $\hat{x}(0)=[1,\ 2,\ 0]^{\top}$. We compare the proposed S-ADRC with the conventional ADRC whose derivative gains are constant $u_{\mathrm{ADRC}}=-p\hat{x}_1-\bar{d}\hat{x}_2-\hat{x}_3$. For the case study of this example, an ADRC with $\bar{d}=4.7$ has the best response in terms of the cost function \eqref{eq:III-8}, which we call it the most efficient ADRC. Fig. \ref{fig:1} shows the response of the system under the proposed S-ADRC. As shown, although the S-ADRC does not have any knowledge of the system model $f$, its response is closer than the most efficient ADRC to the optimal solution obtained via a full knowledge of the system model. As the system is linear in this example, its optimal controller is a linear quadratic regulator (LQR).
\end{example}

\begin{figure}[!h]
\centering\includegraphics[width=1\linewidth]{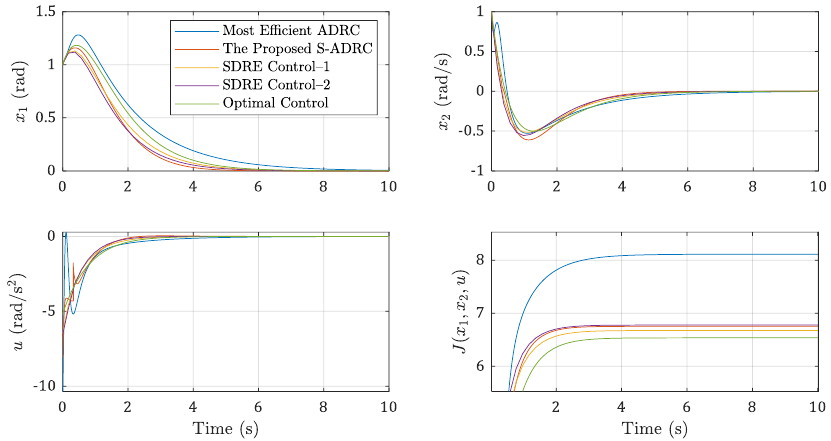}
\caption{Simulation results for Example \ref{ex:2}. SDRE control--1 and 2 correspond to the SDC matrices \eqref{ex2-1a} and \eqref{ex2-1b}, respectively.}
\label{fig:2}
\end{figure}

\begin{example}
\label{ex:2}
In this example, we consider an inverted pendulum with a time-varying point mass whose dynamics can be expressed as \eqref{eq:1} with $f(x_1,t)=(1+\exp(-t))\sin(x_1)$ and $g=1$. The cost function weights are considered as $q_1=q_2=r=1$. The LESO parameters and the system initial conditions are the same as those in Example \ref{ex:1}. The controller parameters are selected to be $p=2$, $\rho=\sqrt{3}$, and $\beta=1+2(p+\rho)$ (corresponding to weight factors of $q_1/r=q_2/r=1$). We compare the response of the model-free S-ADRC with two model-based SDRE controllers \cite{ccimen2008state,mracek1998control} with different SDC matrices obtained as follows:
\begin{subequations}
\label{ex2-1}
\begin{align}
\label{ex2-1a}
\text{SDC matrix--1:}\ A(x)&=\begin{bmatrix}
0 & 1 \\
\frac{(1+\exp(-t))\sin(x_1)}{x_1} & 0
\end{bmatrix},\\
\label{ex2-1b}
\text{SDC matrix--2:}\ A(x)&=\begin{bmatrix}
0 & 1 \\
\sqrt{3} & \frac{(1+\exp(-t))\sin(x_1)-\sqrt{3}x_1}{x_2}
\end{bmatrix},
\end{align}
\end{subequations}
where the SDC matrix--2 is similar to \eqref{eq:th:opt-5}. Similarly to Example \ref{ex:1}, the results are also compared with the most efficient ADRC, $u_{\mathrm{ADRC}}=-p\hat{x}_1-\bar{d}\hat{x}_2-\hat{x}_3$, which has the lowest cost for $\bar{d}=17.4$. As shown in Fig. \ref{fig:2}, the cost of the proposed S-ADRC is lower than that of the most efficient ADRC. Moreover, the S-ADRC's cost value is close to those of the SDRE controllers and the optimal solution, while they have access to the full knowledge of the system model and state variables. The optimal controller is obtained numerically using the variation of extremals algorithm (see \cite[\S 6.3]{kirk2004optimal}).
\end{example}

\begin{figure}[!h]
\centering\includegraphics[width=1\linewidth]{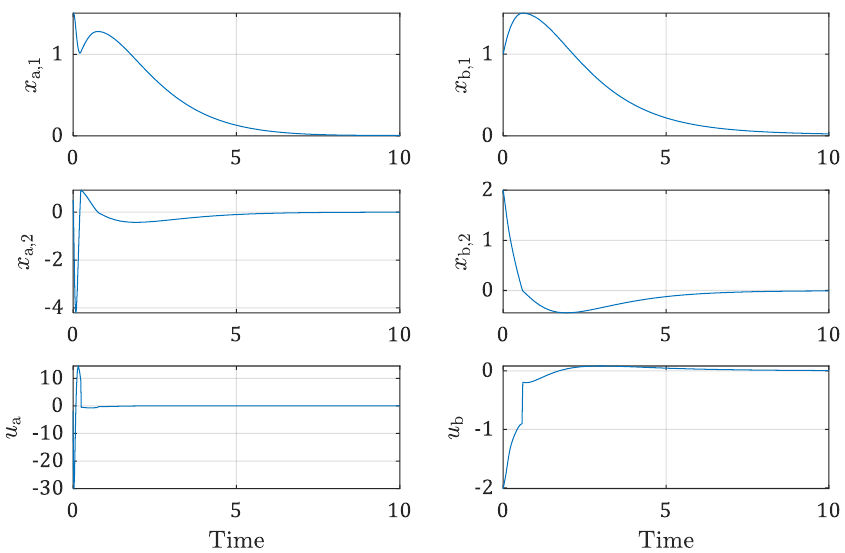}
\caption{Simulation results for Example \ref{ex:3}.}
\label{fig:4}
\end{figure}

\begin{example}
\label{ex:3}
This example discusses the extension of the proposed S-ADRC to multi-variable systems. We consider a two-mass system connected by an unknown nonlinear spring and damper, illustrated in Fig. \ref{fig:3}, as a coupled multi-variable case study. The dynamics of the system can be expressed as follows:
\begin{equation}
\label{eqex:1}
\left\{\begin{array}{lcl}
\dot{x}_{\mathrm{a},1}&=&x_{\mathrm{a},2}, \\
\dot{x}_{\mathrm{b},1}&=&x_{\mathrm{b},2}, \\
\dot{x}_{\mathrm{a},2}&=&f_{\mathrm{a}}(x_{\mathrm{a},1},x_{\mathrm{a},2},x_{\mathrm{b},1},x_{\mathrm{b},2})+\frac{1}{m_{\mathrm{a}}}u_{\mathrm{a}}, \\
\dot{x}_{\mathrm{b},2}&=&f_{\mathrm{b}}(x_{\mathrm{a},1},x_{\mathrm{a},2},x_{\mathrm{b},1},x_{\mathrm{b},2})+\frac{1}{m_{\mathrm{b}}}u_{\mathrm{b}}, \\
y_{\mathrm{a}}&=&x_{\mathrm{a},1}, \\
y_{\mathrm{b}}&=&x_{\mathrm{b},1},
\end{array}\right.
\end{equation}
where $m_{\mathrm{a}}$ and $m_{\mathrm{b}}$ are known, and the unknown dynamics of the system has the following formula:
\begin{subequations}
\label{eq:ex:3-2}
\begin{align}
\begin{split}
f_{\mathrm{a}}=&k_1(x_{\mathrm{b},1}-x_{\mathrm{a},1})+k_2(x_{\mathrm{b},1}-x_{\mathrm{a},1})^3\\&+c_1(x_{\mathrm{b},2}-x_{\mathrm{a},2})+c_2(x_{\mathrm{b},2}-x_{\mathrm{a},2})^3 \end{split},\\
f_{\mathrm{b}}=&-f_{\mathrm{a}}.
\end{align}
\end{subequations}

\begin{figure}[!h]
\centering\includegraphics[width=0.75\linewidth]{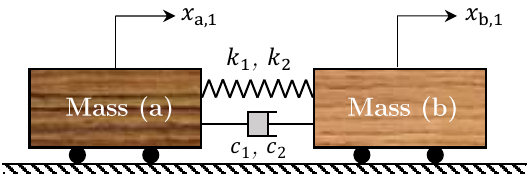}
\caption{A schematic view of the two-mass system of Example \ref{ex:3}.}
\label{fig:3}
\end{figure}
For this system, we extend the single-variable S-ADRC \eqref{eq:III-3} to the following multi-variable form:
\begin{subequations}
\label{eq:ex:3-3}
\begin{align}
u_{\mathrm{a}}&=-\frac{1}{m_{\mathrm{a}}}(p_{\mathrm{a}}\hat{x}_{\mathrm{a},1}+d_{\mathrm{a}}(\hat{x}_\mathrm{a})\hat{x}_{\mathrm{a},2}+\hat{x}_{\mathrm{a},3}),\\
u_{\mathrm{b}}&=-\frac{1}{m_{\mathrm{b}}}(p_{\mathrm{b}}\hat{x}_{\mathrm{b},1}+d_{\mathrm{b}}(\hat{x}_\mathrm{b})\hat{x}_{\mathrm{b},2}+\hat{x}_{\mathrm{b},3}),
\end{align}
\end{subequations}
where two LESOs for each of the masses (a) and (b) work in parallel. We consider the parameters and initial conditions of both LESOs the same as those in Example \ref{ex:1}. The control parameters are selected as $p=1$, $\rho=0.5$, $\beta=4$, and $\delta=10^6$. Fig. \ref{fig:4} shows the simulation results for $k_1=0.5$, $k_2=0.01$, $c_1=0.05$, $c_2=0.002$, $m_{\mathrm{a}}=0.33$, $m_{\mathrm{b}}=0.4$, $x_{\mathrm{a}}(0)=[1.5,\ 0.5]^\top$, and $x_{\mathrm{b}}(0)=[1,\ 2]^\top$. It can be verified that the controller successfully leads the system to the desired state values, i.e., zero position and velocity.  

\end{example}

\section{Conclusions}    
\label{sec:VI}

A suboptimal active disturbance rejection controller (S-ADRC) was developed for second-order nonlinear systems without using their dynamics model. It was proved that the proposed control scheme globally converges to the vicinity of an optimal solution. Therefore, the controller benefits from the model-free and globally stable features of an ADRC, while also acting as an SDRE-based controller in approaching an optimal solution. It was demonstrated through numerical simulations that the cost of the proposed S-ADRC is close to that of the model-based optimal and suboptimal controllers. Also, it was shown that the S-ADRC is more efficient than any conventional ADRC with constant gains. Although the convergence and suboptimality of the S-ADRC were proved for single-variable, second-order nonlinear systems, the applicability of the controller to multi-variable cases was also illustrated through an example. 

The proposed control scheme can be further extended in future works to higher-order systems with uncertain input gain.




\bibliographystyle{IEEEtran}
\bibliography{biblo}

\end{document}